\documentclass[12pt]{article}

\usepackage{hyperref}
\usepackage{amsmath}
\usepackage{amssymb}
\usepackage{amsfonts}
\usepackage{amsopn}
\usepackage{amsthm}
\usepackage[all]{xy}

\swapnumbers
\theoremstyle{plain}
\newtheorem{thm}{Theorem}[section]
\newtheorem{lem}[thm]{Lemma}

\newtheorem{prop}[thm]{Proposition}
\newtheorem*{mainthm}{Theorem}

\theoremstyle{definition}
\newtheorem{ntt}[thm]{}

\newcommand{\pl}{\mathbb{P}}   
\newcommand{\zz}{\mathbb{Z}}   
\newcommand{\D}{\mathrm{D}} 
\newcommand{\E}{\mathrm{E}} 
\newcommand{\F}{\mathrm{F}_4} 
\newcommand{\HH}{\mathrm{H}^1} 

\newcommand{\M}{\mathcal{M}}     
\newcommand{\OP}{\mathbb{OP}^2} 
\newcommand{\X}{X} 
\newcommand{\B}{\mathcal B} 
\newcommand{\RatCor}{RatCor}    

\DeclareMathOperator{\Spec}{\mathrm{Spec}}  
\DeclareMathOperator{\Spin}{\mathrm{Spin}}  
\DeclareMathOperator{\Gal}{\mathrm{Gal}}   
\DeclareMathOperator{\CH}{\mathrm{CH}}      
\DeclareMathOperator{\CHO}{\widetilde{\mathrm{CH}}_0}
\DeclareMathOperator{\Aut}{\mathrm{Aut}}    
\DeclareMathOperator{\Sim}{\mathrm{Sim}}    
\DeclareMathOperator{\Mor}{\mathrm{Mor}}    
\DeclareMathOperator{\Ker}{\mathrm{Ker}}    
\DeclareMathOperator{\im}{\mathrm{Im}}      
\DeclareMathOperator{\SB}{\mathrm{SB}}      
\DeclareMathOperator{\PGL}{\mathrm{PGL}}      
\DeclareMathOperator{\GO}{\mathrm{GO}}  
\DeclareMathOperator{\Gr}{\mathrm{Gr}} 
\DeclareMathOperator{\Stab}{\mathrm{Stab}}  
\DeclareMathOperator{\diag}{\mathrm{diag}} 

\title{Zero cycles on a twisted Cayley plane}

\author{V.~Petrov, N.~Semenov, K.~Zainoulline}
\date{}

\begin{document}

\maketitle

\begin{abstract}
Let $k$ be a field of characteristic not 2 and 3.
Let $G$ be an exceptional simple algebraic group over $k$
of type $\F$, $^1{\E_6}$ or $\E_7$ with trivial Tits algebras.
Let $X$ be a projective $G$-homogeneous variety.
If $G$ is of type $\E_7$ we assume in addition
that the respective
parabolic subgroup is of type $P_7$.
The main result of the paper says that 
the degree map on the group of zero cycles of $X$ 
is injective.
\end{abstract}

\section{Introduction}

Let $k$ be a field and $G$ a simple algebraic group over $k$.
Consider a projective $G$-homogeneous variety $X$ over $k$.
Any such variety over the separable closure $k_s$
of $k$ becomes isomorphic to the quotient $G_s/P$,
where $P$ is a parabolic subgroup of the split group $G_s=G\times_k k_s$.
It is known that the conjugacy classes of parabolic subgroups
of $G_s$ are in one-to-one correspondence with the subsets
of the vertices $\Pi$ of the Dynkin diagram of $G_s$:
we say a parabolic subgroup is of type $\theta \subset \Pi$ and
denote it by $P_\theta$ if it is conjugate to a standard parabolic subgroup
generated by the Borel subgroup and all unipotent subgroups
corresponding to roots in the span of $\Pi$ with no
$\theta$ terms (see \cite[42.3.1]{TW02}).

In the present paper we assume $k$ has characteristic not $2$ and $3$.
$G$ is an exceptional simple algebraic group over $k$ 
of type $\F$, $^1{\E_6}$ or $\E_7$ with trivial Tits algebras and
$X$ is a projective $G$-homogeneous variety over $k$.
The goal of the paper is to compute the group of zero-cycles 
$\CH_0(X)$ which is an important geometric invariant of a variety.
Namely, we prove

\begin{mainthm}
Let $k$ be a field of characteristic not $2$ and $3$.
Let $G$ be an exceptional simple algebraic group over $k$ 
of type $\F$, $^1{\E_6}$ or $\E_7$ with trivial Tits algebras and
$X$ a projective $G$-homogeneous variety over $k$.
If $G$ is of type $\E_7$ we assume in addition that
$X$ corresponds to the parabolic subgroup of type $P_7$.
Then the degree map $\CH_0(X) \to \zz$ is injective.
\end{mainthm}

The history of the question starts with the work of I.~Panin \cite{Pa84},
where he proved the injectivity of the degree map for Severi-Brauer varieties.
For quadrics this was proved by R.~Swan in \cite{Sw89}. 
The case of involution varieties was considered by A.~Merkurjev in \cite{Me95}. 
For varieties of type $\F$ it was announced by M.~Rost.

Our work was mostly motivated by the paper of D.~Krashen \cite{Kr05}, where
he reformulated the question above in terms of $R$-triviality 
of certain symmetric powers and proved the injectivity
for a wide class of generalized Severi-Brauer varieties 
and some involutive varieties,
hence, generalizing the previously known results of Panin and Merkurjev.
Another motivating point was the result of V.~Popov \cite{Po05},
giving a full classification
of generically $n$-transitive actions 
of a split linear algebraic group $G$ on
a projective homogeneous variety $G/P$. 
For instance, the case of a Cayley plane $X=G/P_1$, where
$G$ is split of type $\E_6$ (see \cite{IM05}), 
provides an example of such an action for $n=3$.
As a consequence, one can identify the open orbit $X^0$ 
of this action with a
homogeneous variety $G/(T\cdot \Spin_8)$, where 
$T$ is the torus which is complementary to $\Spin_8$.
Then the results of Krashen reduce the question of injectivity to the 
question of $R$-triviality of the twisted form of $X^0/S_3$.

Apart from the main result concerning exceptional varieties 
we give shortened proofs for injectivity
of the degree map for quadrics and Severi-Brauer varieties as well.

As mentioned to us by V.~Chernousov and M.~Rost,
it is possible to prove the same results
using Rost invariant and the Chain Lemma. Our proof doesn't use these tools but
only geometry and several classical results 
concerning exceptional groups.

The paper is organized as follows.
In the first section we provide several facts about zero-cycles
and symmetric powers. Then we prove the theorem for twisted forms
of a Cayley plane (here prime 3 plays the crucial role).
In the next section we prove injectivity for the twisted
form of a homogeneous variety of type $\E_7$ (this deals
with prime 2). In the last section, 
we combine these two results together with some facts
about rational correspondences and finish the proof of the main theorem.

\section{Zero cycles and symmetric powers}
In this section we collect some facts from \cite{Kr05} on the interrelation between zero cycles on a projective variety $X$ and classes of $R$-equivalence on symmetric powers of $X$.

\begin{ntt}\label{descent}
We systematically use Galois descent language, i.e., identify a 
(quasi-projective) variety $X$ over $k$ with the variety $X_s=X\times_k k_s$ 
over the separable closure $k_s$ equipped with an action 
(by semiautomorphisms of $X_s$) of the absolute Galois group 
$\Gamma=\Gal(k_s/k)$. This means that for any $\sigma\in\Gamma$ we are 
given an automorphism $\varphi_\sigma$ of $X_s$ over $k$ such that
the diagram
$$
\xymatrix{
X_s\ar[d]\ar[r]^{\varphi_\sigma}&X_s\ar[d]\\
\Spec k_s\ar[r]^{\sigma^\sharp}&\Spec k_s
}
$$
is commutative, and  $\varphi_{\sigma\tau}=\varphi_\sigma\varphi_\tau$.

The set of $k$-rational points of $X$ is precisely the set of 
$k_s$-rational points of $X_s$ stable under 
the action of $\Gamma$.
\end{ntt}

\begin{ntt}
Let $X$ be a variety over $k$. 
Two rational points $p,q\in X(k)$ are called 
\emph{elementary linked} if there exists a rational morphism 
$\varphi\colon\pl^1_k\dasharrow X$ such that $p,q\in\im(\varphi(k))$. 
The $R$-equivalence is the equivalence relation generated by this relation. 
A variety $X$ is called \emph{$R$-trivial} if there exists exactly one class 
of $R$-equivalence on $X$, and \emph{algebraically $R$-trivial} if 
$X_K=X\times_kK$ is $R$-trivial for any finite field extension 
$K/k$.

The $n$-th symmetric power of $X$ is by the definition the variety 
$S^nX=X^n/S_n$, 
where $S_n$ is the symmetric group acting on $X^n$ via permutations.

Let $p$ be a prime number. A field $k$ is called \emph{prime-to-$p$ closed} 
if there is no proper, finite field extension $K/k$ of degree prime to $p$. 
For any field $k$ we denote by $k_p$ its \emph{prime-to-$p$ closure} that is the algebraic extension of $k$ which is prime-to-$p$ closed.

Let $X$ be a projective variety over $k$. By $\CHO(X)$ we denote
the kernel of the degree map:
$$\CHO(X)=\Ker(\deg\colon\CH_0(X)\to\zz).$$
\end{ntt}

The following results will play a crucial role in the sequel:

\begin{prop}\label{propPrime}(\cite[Lemma~1.3]{Kr05}) If $\CHO(X_{k_p})=0$ for each prime $p$ then $\CHO(X)=0$.
\end{prop}

\begin{prop}\label{propKr}(\cite[Theorem~1.4]{Kr05}) Suppose that $k$ is 
prime-to-$p$ closed and the following conditions are satisfied:
\begin{enumerate}
\item $S^{p^n}X$ is algebraically $R$-trivial for some $n\ge 0$,
\item For any field $K/k$ such that $X(K)\ne\emptyset$ the variety $X_K$ is $R$-trivial.
\end{enumerate}
Then $\CHO(X)=0$.
\end{prop}

\begin{ntt}\label{SBr}
As an easy application of these results 
we sketch a proof of the well-known result 
(first appearing in \cite[Theorem~2.3.7]{Pa84}) 
that $\CHO(\SB(A))=0$, where $A$ 
is a central simple algebra over $k$
(a more general case of flag varieties is settled in \cite{Kr05}). 
For simplicity we assume $\deg A=p$ is prime.

By Proposition~\ref{propPrime} we may assume the base field $k$ 
is prime-to-$p$ closed 
(for a prime $q$ different from $p$ the algebra $A$ splits over $k_q$). 
According to Proposition~\ref{propKr} 
it suffices to show that $S^p\SB(A_K)$ is $R$-trivial 
for every finite extension $K/k$ 
(the second hypothesis of \ref{propKr} 
holds for any twisted flag variety). 
Changing the base we may assume $K=k$. 
If $A$ is split, the assertion is trivial; so we may assume $A$ is not split.

According to our conventions (see \ref{descent}) 
the variety $\SB(A)$ is the 
variety of all parabolic subgroups $P$ of type $P_1$ in 
the group $\PGL_1(A\otimes_kk_s)$ 
with the action of $\Gamma$ coming from its 
action on $k_s$. 
Therefore, $S^pX$ is the variety 
of all unordered $p$-tuples $[P^{(1)},\ldots,P^{(p)}]$ 
of parabolic subgroups of type $P_1$ of 
$\PGL_1(A\otimes_kk_s)$. 
Let $U$ be an open subset of $S^pX$ 
defined by the condition that 
the intersection $P^{(1)}\cap\ldots\cap P^{(p)}$ 
is a maximal torus in $\PGL_1(A\otimes_kk_s)$. 
Every maximal torus $T$ in $\PGL_1(A\otimes_kk_s)$ 
is contained in precisely $p$ parabolic subgroups of type $P_1$, 
whose intersection is $T$. Therefore, 
$U$ is isomorphic to the variety of all maximal tori in $\PGL_1(A)$. 
This variety is known to be rational 
(and therefore $R$-trivial since it is homogeneous). 
Moreover, one can check that if $A$ is not split 
then the embedding $U\to S^pX$ is surjective on $k$-points. 
So $S^pX$ is $R$-trivial, and we are done.
\end{ntt}

\begin{ntt}\label{quadr}
The same method can be applied to prove that 
$\CHO(Q)=0$ for a nonsingular projective quadric 
$Q$ over a field of characteristic not $2$ 
(the result of Swan \cite{Sw89}).

As above, we may assume that $p=2$ and $Q$ is anisotropic.
It suffices to prove  
that $S^2Q$ is $R$-trivial. 
Let $q$ be the corresponding quadratic form on a 
vector space $V$. 
The quadric $Q$ can be viewed as the variety of lines 
$\langle v\rangle$, 
where $v\in V\otimes_kk_s$ satisfies 
$q(v)=0$, with the obvious action of $\Gamma$. 
So $S^2Q$ can be identified with 
the variety of pairs 
$[\langle v_1\rangle,\langle v_2\rangle]$ of lines of this kind, 
with induced action of $\Gamma$. 
Consider the open subset $U$ defined by the condition 
$b_q(v_1,v_2)\ne 0$ ($b_q$ stands for the polarization of $q$). 
Clearly, the embedding $U\to S^2Q$ is surjective on $k$-points 
(otherwise the subspace $\langle v_1,v_2\rangle$ 
defines a totally isotropic subspace over $k$). 
So it is enough to check that $U$ is $R$-trivial.

Consider the open subvariety $W$ of $\Gr(2,V)$ 
consisting of planes
$H\subset V\otimes_kk_s$ such that $q|_H$ is nonsingular. 
For every such a plane there exists up to scalar factors 
exactly one hyperbolic basis $\{v_1,v_2\}$ over $k_s$. 
Therefore, the map from $U$ to $W$ sending 
$[\langle v_1\rangle,\langle v_2\rangle]$ to $\langle v_1,v_2\rangle$ 
is an isomorphism. 
But any open subvariety of $\Gr(2,V)$ is $R$-trivial, and we are done.
\end{ntt}

We shall use in the sequel the following observation.
\begin{lem}\label{lemSur}
Let $H\subset K\subset G$ be algebraic groups over $k$. 
Suppose that the map $\HH(k,H)\to\HH(k,K)$ is surjective. 
Then the morphism $G/H\to G/K$ is surjective on $k$-points.
\end{lem}
\begin{proof}
An element $x$ of $G/K(k)$ is presented by an element $g\in G(k_s)$ 
satisfying the
condition $\gamma(\sigma)=g^{-1}\cdot{}^\sigma g$ lies in $K(k_s)$ 
for all $\sigma\in\Gamma$. But $\gamma$ is clearly a $1$-cocycle 
with coefficients in $K$. Therefore by assumption, there exists some $h\in K$ such that $h^{-1}\gamma(\sigma)\cdot{}^\sigma h=(gh)^{-1}\cdot{}^\sigma(gh)$ is a $1$-cocycle with coefficients in $H$. But then $gh$ presents an element of $G/H(k)$ which goes to $x$ under the morphism $G/H\to G/K$, and the lemma is proved.
\end{proof}

\section{Twisted forms of a Cayley plane}
In the present section we prove the injectivity of the degree map 
in the case when $X$ is a twisted form of a Cayley plane.

\begin{ntt}
Let $J$ denote a simple exceptional $27$-dimensional Jordan algebra over $k$, 
and $N_J$ its norm (which is a cubic form on $J$). 
An invertible linear map $f\colon J\to J$ is called a \emph{similitude} if there exists some $\alpha\in k^*$ (called the \emph{multiplier} of $f$) such that $N_J(f(v))=\alpha N_J(v)$ for all $v\in J$. The group
$G=\Sim(J)$ of all similitudes is a reductive group of type $^1{\E_6}$, and 
every group of type $^1{\E_6}$ with trivial Tits algebras can be obtained in this way up to isogeny (see \cite[Theorem~1.4]{Ga01}).
\end{ntt}

\begin{ntt} The \emph{(twisted) Cayley plane} $\OP(J)$ 
is the variety obtained by Galois descent from the variety of all parabolic subgroups of type $P_1$ in $\Sim(J\otimes_k k_s)$ (see \cite{IM05}).
This variety can be identified with the variety of all lines $\langle e\rangle$ spanned by elements
$e\in J_s=J\otimes_k k_s$ satisfying the condition $e\times e=0$ (see \cite[Theorem~7.2]{Ga01}). 
\end{ntt}

The goal of the present section is to prove
\begin{thm}\label{cayl}
$\CHO(\OP(J))=0$.
\end{thm}

We start the proof with the following easy reduction.

\begin{ntt} By Proposition~\ref{propKr} it is enough to prove that 
$(S^p \OP(J))\times_k K$ is $R$-trivial for any prime $p$ and any finite 
field extension $K/k_p$.
Making the base change, it suffices to prove it for $K=k_p$.
Moreover, we may assume $J$ is not reduced 
(otherwise $\OP(J)$ is a rational homogeneous variety and, hence, is 
$R$-trivial).

Assume $p\ne 3$, then $\OP(J)(k_p)\ne\emptyset$ (and, hence, is $R$-trivial). 
Indeed, choose any cubic \'etale subalgebra $L$ of $J$ 
(see \cite[Proposition~39.20]{Inv}). It splits over $k_p$ and, 
therefore, $L\otimes_kk_p$ contains a primitive idempotent $e$. 
As an element of $J\otimes_kk_p$ it satisfies the condition $e\times e=0$ 
(see \cite[Lemma~5.2.1(i)]{SV}).
So we may assume $p=3$. 
\end{ntt}

\begin{ntt}
>From now on $p=3$ and the field $k$ is prime-to-$p$ closed.
By definition $S^3(\OP(J))$ is the variety of all unordered triples
$[\langle e_1\rangle,\langle e_2\rangle,\langle e_3\rangle]$, 
where $e_i$ are the elements of $J_s=J\otimes_kk_s$ satisfying the conditions 
$e_i\times e_i=0$, with the natural action of $\Gamma$. 
Denote by $U$ the open subvariety of $\OP(J)$ 
defined by the condition
$$
N_{J_s}(e_1,e_2,e_3)\ne 0,
$$
where $N$ is the polarization of the norm.

The embedding $U\to S^3(\OP(J))$ is surjective on $k$-points. 
Indeed, if
$[\langle e_1\rangle,\langle e_2\rangle,\langle e_3\rangle]$ 
is stable under the action of $\Gamma$ and
$N_{J_s}(e_1,e_2,e_3)=0$, then 
$\langle e_1,e_2,e_3\rangle$ gives by descent a 
$k$-defined subspace $V$ of $J$ such that $N|_V=0$. 
But then $J$ is reduced by \cite[Theorem~5.5.2]{SV}, 
which leads to a contradiction. 
So it is enough to show that $U$ is $R$-trivial.
\end{ntt}

\begin{ntt}
Choose a cubic \'etale subalgebra $L$ in $J$. 
Over the separable closure we have
$$
L\otimes_kk_s=k_s e_1\oplus k_s e_2\oplus k_s e_3
$$
where $e_1,e_2,e_3$ are primitive idempotents. 
Then we have $e_i\times e_i=0$ in $J_s$, $i=1,2,3$, the norm
$N_{J_s}(e_1,e_2,e_3)=N_{L\otimes_kk_s}(e_1,e_2,e_3)$ is trivial
and the triple $[e_1,e_2,e_3]$ is stable under the action of 
$\Gamma$ (since so is $L$). 
Hence, the triple $[\langle e_1\rangle,\langle e_2\rangle,\langle e_3\rangle]$ 
is a $k$-rational point of $U$.

By \cite[Proposition~3.12]{SV68} the group 
$G$ acts transitively on $U$. Therefore, we have
$$
U\simeq G/\Stab_G([\langle e_1\rangle,\langle e_2\rangle,\langle e_3\rangle]).
$$
The stabilizer is defined over $k$, since it is stable under $\Gamma$.
Moreover, it
coincides with $\Stab_G(L)$. Indeed, one inclusion is obvious, 
and the other one follows from the fact that $e_1,e_2,e_3$ are 
the only elements $e$ of $L\otimes_kk_s$ satisfying the condition 
$e\times e=0$, up to scalar factors (see \cite[Theorem~5.5.1]{SV}).
\end{ntt}

\begin{ntt}\label{descrex}
Consider the Springer decomposition $J=L\oplus V$ of $J$ with respect to 
$L$. The pair $(L,V)$ has a natural structure of a 
twisted composition, and there is a monomorphism
$\Aut(L,V)\to\Aut(J)$ sending a pair $(\varphi,t)$ 
(where $\varphi\colon L\to L$,
$t\colon V\to V$) to $\varphi\oplus t\colon J\to J$ 
(see \cite[\S~38.A]{Inv}). Note that
$\Aut(L,V)$ coincides with the stabilizer of $L$ in $\Aut(J)$.
\end{ntt}

\begin{lem} The following sequence of algebraic groups is exact
\begin{align*}
1\longrightarrow\Aut(L,V)\longrightarrow\Stab_G(L)&\longrightarrow
R_{L/k}(\mathbb G_m)\longrightarrow 1,\\
f&\mapsto f(1)
\end{align*}
where $R_{L/k}$ stands for the Weil restriction. 
\end{lem}

\begin{proof}
Exactness at the middle term follows 
from \ref{descrex} and 
the fact that the stabilizer of $1$ 
in $G$ coincides with $\Aut(J)$ (see \cite[Proposition~5.9.4]{SV}).
To prove the exactness at the last term observe that 
a $k_s$-point of $R_{L/k}(\mathbb G_m)$ 
is a triple of scalars $(\alpha_0,\alpha_1,\alpha_2)\in k_s^*\times k_s^*\times k_s^*$. Multiplying $f$ by the scalar transformation with the coefficient
$(\alpha_0\alpha_1\alpha_2)^\frac{1}{3}$ 
(which is an element of $\Stab_G(L)(k_s)$) we may assume that 
$\alpha_0\alpha_1\alpha_2=1$. Choose a \emph{related triple}
$(t_0,t_1,t_2)$ of elements of $\GO^+(\mathbb O_d,N_{\mathbb O_d})$ ($\mathbb O_d$ is the split Cayley algebra) such that 
$\mu(t_i)=\alpha_i$, $i=0,1,2$ (see \cite[Corollary~35.5]{Inv}). 
Now the transformation $f$ of $J$ defined by
$$
\begin{pmatrix}
\varepsilon_0&c_2&\cdot\\
\cdot&\varepsilon_1&c_0\\
c_1&\cdot&\varepsilon_2
\end{pmatrix}
\mapsto
\begin{pmatrix}
\alpha_0\varepsilon_0&t_2(c_2)&\cdot\\
\cdot&\alpha_1\varepsilon_1&t_0(c_0)\\
t_1(c_1)&\cdot&\alpha_2\varepsilon_2
\end{pmatrix}
$$
lies in $\Sim(J)$ by \cite[(7.3)]{Ga01}, stabilizes $L_{k_s}=\diag(k_s,k_s,k_s)$ and sends $1\in J_s$ to $\diag(\alpha_0,\alpha_1,\alpha_2)$.
\end{proof}

\begin{ntt}
Since $\HH(k,L^*)=1$ (by Hilbert '90), 
the map $\HH(k,\Aut(L,V))\to\HH(k,\Stab_G(L))$ is surjective. 
By Lemma~\ref{lemSur} the morphism
$$
G/\Aut(L,V)\to G/\Stab_G(L)\simeq U
$$
is surjective on $k$-points. Therefore, 
it suffices to show that $G/\Aut(L,V)$ is $R$-trivial.
\end{ntt}

\begin{ntt}
Consider the morphism
$$
\psi\colon G/\Aut(L,V)\to G/\Aut(J).
$$
By \cite[Corollary~3.14]{Kr05} it suffices to show that
\begin{enumerate}
\item $\psi$ is surjective on $k$-points;
\item $G/\Aut(J)$ is $R$-trivial;
\item The fibers of $\psi$ (which are isomorphic to $\Aut(J)/\Aut(L,V)$) are unirational and $R$-trivial.
\end{enumerate}
\end{ntt}

\begin{ntt}
In order to prove surjectivity of $\psi$ on $k$-points 
it is enough by Lemma~\ref{lemSur} to prove surjectivity 
of the map $\HH(k,\Aut(L,V))\to\HH(k,\Aut(J))$. The set 
$\HH(k,\Aut(L,V))$ classifies all twisted compositions 
$(L',V')$ which become isomorphic to $(L,V)$ over 
$k_s$ and $\HH(k,\Aut(J))$ classifies all 
(exceptional $27$-dimensional) Jordan algebras $J'$. 
It is easy to verify that the morphism sends $(L',V')$ 
to the Jordan algebra $L'\oplus V'$ and, hence, 
the surjectivity follows from the fact that any Jordan algebra 
admits a Springer decomposition (cf. \cite[Proposition~38.7]{Inv}).
\end{ntt}

\begin{ntt}
Let $W$ be the open subvariety of $J$ 
consisting of elements $v$ with $N_J(v)\ne 0$. 
Then $G$ acts transitively on $W$ (see \cite[Proposition~5.9.3]{SV}) 
and the stabilizer of the point $1$ coincides with $\Aut(J)$.
So $G/\Aut(J)\simeq W$ is clearly $R$-trivial.
\end{ntt}

\begin{ntt}
Consider the variety $Y$ of all \'etale cubic subalgebras of $J$. 
By \cite[Proposition~39.20(1)]{Inv} 
there is a map from an open subvariety $J_0$ 
of regular elements in $J$ to $Y$ (sending $a$ to $k[a]$), 
surjective on $k$-points. Therefore $Y$ is unirational and $R$-trivial
$24$-dimensional irreducible variety.

The group $\Aut(J)$ acts on $Y$ naturally. 
Let $L'$ be any $k$-point of $Y$. 
The stabilizer of $L'$ in $\Aut(J)$ obviously equals to 
$\Aut(L',V')$ ($J=L'\oplus V'$ is the Springer decomposition). 
So the orbit of $L'$ is isomorphic to $\Aut(J)/\Aut(L',V')$ and, 
in particular, has dimension $24$. Therefore, it is open and, since $L'$ is arbitrary, the action is transitive. 
So we have $\Aut(J)/\Aut(L,V)\simeq Y$ is unirational and $R$-trivial 
and the proof of \ref{cayl} is completed.
\end{ntt}

\section{Case $\E_7/P_7$}
In the present Section we prove the injectivity of the degree map
for twisted forms of 
a projective homogeneous variety corresponding
to an exceptional group of type $\E_7$ and
a parabolic subgroup of type $P_7$.

\begin{ntt}
Let $\B$ denote a $56$-dimensional Brown algebra over $k$. It defines naturally up to a scalar factor a skew-symmetric form $b$ on $\B$ and a trilinear map $t$ from $\B\times\B\times\B$ to $\B$ such that $(\B,t,b)$ is a Freudenthal triple system (see \cite[Definition~3.1]{Ga01} and \cite[\S~4]{Ga01}). An invertible linear map $f\colon \B\to\B$ is called a \emph{similitude} if there exists some $\alpha\in k^*$ (called the \emph{multiplier} of $f$) such that $b(f(u),f(v))=\alpha b(u,v)$ and
$t(f(u),f(v),f(w))=\alpha f(t(u,v,w))$ for all $u,v,w\in\B$. The group
$G=\Sim(\B)$ of all similitudes is a reductive group of type $\E_7$ 
and every group of 
type $\E_7$ with trivial Tits algebras can be obtained in this way up to 
isogeny (cf. \cite[Theorem~4.16]{Ga01}).

An element $e$ is called \emph{singular} (or \emph{strictly regular} 
following \cite{Fe72}) if $t(e,e,\B)\subseteq\langle e\rangle$. 
In this case $t(e,e,v)=2b(v,e)e$ for every $v\in V$. 
An equivalent definition is that $t(e,e,e)=0$  and $e\in t(e,e,\B)$ 
(see \cite[Lemma~3.1]{Fe72}). $\B$ is called \emph{reduced} if it 
contains singular elements. There do exist anisotropic groups of type $\E_7$
with trivial Tits algebras over certain fields (see \cite{Ti90}).
\end{ntt}

\begin{ntt}
Let $\X(\B)$ be the variety obtained by the Galois descent 
from the variety of all parabolic subgroups of type $P_7$  
in $\Sim(\B\otimes_kk_s)$ 
(the action of $\Gamma=\Gal(k_s/k)$ comes from the action on $k_s$). 
This variety can be identified with the variety of all lines 
$\langle e\rangle$ spanned by singular elements
$e\in\B\otimes_k k_s$ (see \cite[Theorem~7.6]{Ga01}).
\end{ntt}

The goal of this section is to prove

\begin{thm}\label{E7case}
$\CHO(\X(\B))=0$.
\end{thm}

We start with the similar reduction as in the case of $\E_6$

\begin{ntt}\label{qwer}
Assume first that $G$ has Tits index $\E_{7,1}^{66}$ 
(see \cite[Table~II]{Ti66}). 
Its anisotropic kernel is of type
$\D_6$ and, since $G$ has trivial Tits algebras, 
the anisotropic kernel corresponds to a 
$12$-dimensional nondegenerate quadratic form $q$ with 
split simple factors of its Clifford algebra. 
A straightforward computation (see \cite[Thm.~7.4]{Br05}) shows that
$$\M(\X(\B))\simeq\M(Q)\oplus\M(Y)(6)\oplus\M(Q)(17),$$
where $Q$ is the projective quadric corresponding to $q$, $Y$ is a twisted
form of the maximal orthogonal grassmanian of a split $12$-dimensional
quadric, and $\M$ denotes Chow motive. 
Therefore, $\CHO(\X(\B))=\CHO(Q)=0$, the last equality due to Swan.
\end{ntt}

\begin{ntt} By Proposition~\ref{propKr} it is enough to prove that 
$(S^p \X(\B))\times_k K$ is $R$-trivial for any prime $p$ and any finite 
field extension $K/k_p$.
Making the base change, it suffices to prove it for $K=k_p$.
Moreover, we may assume $\B$ is not reduced 
(otherwise $\X(\B)$ is rational and, hence, $R$-trivial).

Assume $p\ne 2$, then $\B\otimes k_p$ is reduced by 
\cite[Corollary~3.4]{Fe72} and, therefore, $\X(\B)(k_p)\ne\emptyset$. 
So we may assume $p=2$. 

>From now on $p=2$ and $k=k_p$.
Since $\B$ is not reduced, the group $G$ has Tits index either 
$\E_{7,0}^{133}$ or $\E_{7,1}^{66}$ 
(see \cite[6.5.5]{Ti71} and \cite[Table~II]{Ti66}).
By \ref{qwer} we may assume $G$ is anisotropic (has index $\E_{7,0}^{133}$).
\end{ntt}

\begin{ntt} By definition $S^2(\X(\B))$ is the variety of all unordered pairs
$[\langle e_1\rangle,\langle e_2\rangle]$, 
where $e_i$ are singular elements of $\B\otimes_kk_s$, 
with the natural action of $\Gamma$. 
Denote by $U$ the open subvariety of $\X(\B)$ 
defined by the condition $b(e_1,e_2)\ne 0$.
\end{ntt}

\begin{lem}
The embedding $U\to S^2(\X(\B))$ is surjective on $k$-points. 
\end{lem}

\begin{proof}
Consider the diagonal action of $G$ on $\X(\B)\times \X(\B)$
(we may assume in this proof that $G$ is simple).
Over $k_s$ this action has four orbits: 
the minimal orbit which is
the diagonal and, hence, is isomorphic to $G_s/P_7$,
the open dense orbit which is isomorphic
to the quotient $G_s/L(P_7)$, where $L(P_7)$ denotes the Levi part of $P_7$, 
and two locally closed orbits.
Indeed, there is a one-to-one
correspondence between the orbits of the $G_s$-action 
and double coset classes $P_7\backslash G_s/P_7$ given by mutually inverse maps
$G_s\cdot (x,y)\mapsto P_7x^{-1}yP_7$ and $P_7wP_7\mapsto G_s\cdot(1,w)$.
Observe that the minimal orbit corresponds to the class of the identity
and the open dense orbit to the class of the longest element
$w_0$ of the Weyl group of $G_s$.

Consider the diagonal action of $G$ on $S^2(\X(\B))$. Over $k_s$
the subset $U$ is the open dense orbit in $S^2(\X(\B))$ (see \ref{lll}).
Assume that there exists a $k$-rational point on $S^2(\X(\B)) \setminus U$.
Then the stabilizer $H$ of this point
is a subgroup of $G$ defined over $k$.
Observe that over $k_s$ the connected component of the identity $H^0$ 
is the stabilizer of one of the non-open orbits for
the action of $G$ on $\X(\B)\times \X(\B)$ considered above, i.e.,
can be identified with the intersection of
two parabolic subgroups $H^0_s=P_7\cap wP_7w^{-1}$, where $w$ is
the double coset representative corresponding to the orbit.
By \cite[Expos\'e~XXVI, Theorem~4.3.2]{DG} 
$H^0_s$ is reductive iff $H^0_s$ is the Levi subgroup of $P_7$,
i.e., iff $P_7wP_7=P_7w_0P_7$.
Therefore, $H^0_s$ is non-reductive and so is $H$. 
The latter implies that $G$ 
must have a unipotent element over $k$.
But according to \cite[p.~265]{Ti86}, if $G$ is anisotropic and $\mathrm{char}\,k\ne 2,3$, 
then this is impossible, a contradiction.
\end{proof}

According to the lemma it suffices to show that $U$ is $R$-trivial.

\begin{ntt}\label{lll}
The Brown algebra $\B\otimes_kk_s$ is split, 
that is isomorphic to the Brown algebra of matrices of the form
$$
\begin{pmatrix}
F&J_d\\
J_d&F
\end{pmatrix},
$$
where $J_d$ is the split Jordan algebra. Set
$$
e_1=\begin{pmatrix}1&0\\0&0\end{pmatrix},\quad e_2=\begin{pmatrix}0&0\\0&1\end{pmatrix}.
$$
The pair $[\langle e_1\rangle,\langle e_2\rangle]$ 
is stable under an arbitrary semiautomorphism of $\B\otimes_kk_s$ 
(see \cite[Proof of Theorem~2.9]{Ga01}) and, in particular, 
under the action of $\Gamma$. Therefore, 
$[\langle e_1\rangle,\langle e_2\rangle]$ is a $k$-rational point of $U$. 
Moreover, $\langle e_1,e_2\rangle$ defines by descent the ($k$-defined) 
\'etale quadratic subalgebra $L$ of $\B$.

By \cite[Proposition~7.6]{Fe72} 
$G$ acts transitively on $U$. Therefore,
$$
U\simeq G/\Stab_G([\langle e_1\rangle,\langle e_2\rangle]).
$$
This stabilizer clearly coincides with $\Stab_G(L)$ (one inclusion is obvious, and the other one follows from the fact that $e_1,e_2$ are the only singular elements of $L\otimes_kk_s$ up to scalar factors).
\end{ntt}

\begin{lem}
There is an exact sequence of algebraic groups
\begin{align*}
1\longrightarrow\Aut(\B)\longrightarrow\Stab_G(L)&\longrightarrow
R_{L/k}(\mathbb G_m)\longrightarrow 1.\\
f&\mapsto f(1)
\end{align*}
\end{lem}

\begin{proof}
This follows from the fact that 
the stabilizer of $1$ in $G$ coincides with $\Aut(\B)$.
Indeed, 
we have an obvious injection $\Aut(\B)\to\Stab_G(1)$. To prove
the surjectivity we can assume that $k$ is separably closed. Let
$f$ be an element of $G$ preserving $1$. Since a
decomposition into a sum of two nonorthogonal singular elements is
unique by \cite[Lemma~3.6]{Fe72} and $1=e_1+e_2$, 
the element $f$ must preserve the pair $[e_1,e_2]$. By
\cite[Lemma 7.5]{Fe72} $f$ has a form $\eta^{\lambda}_{\pi}$, where $\eta$ is
a similitude of $J$ with a multiplier $\rho$, $\pi$ is a permutation
on $\{1,2\}$, $\lambda\in k^*$, and $\eta^{\lambda}_{\pi}$ acts on $\B$ 
by formulae \cite[(15)]{Fe72}. 
Now it follows that $\lambda^{-1}\rho^{-1}=1$ and
$\lambda^2\rho=1$ and therefore $\lambda=\rho=1$. So $f$ is an
automorphism of $\B$, as claimed.

The surjectivity of the last map also follows from \cite[Lemma~7.5]{Fe72}.
\end{proof}

\begin{ntt}
Since $\HH(k,L^*)=1$,
the map $\HH(k,\Aut(\B))\to\HH(k,\Stab_G(L))$ is surjective. 
By Lemma~\ref{lemSur} the morphism
$$
G/\Aut(\B)\to G/\Stab_G(L)\simeq U
$$
is surjective on $k$-points. 
Therefore, it suffices to show that $G/\Aut(\B)$ is $R$-trivial.
\end{ntt}

\begin{ntt}
Let $W$ be the open subvariety of $\B$ consisting of elements $v$ such that
$b(v,t(v,v,v))\ne 0$. Then $G$ acts transitively on $W$ 
(it follows easily from \cite[Theorem~7.10]{Fe72} or \cite[p.~140]{SK77}) 
and the stabilizer of the point $1$ coincides with $\Aut(\B)$. 
So $G/\Aut(\B)\simeq W$ is clearly $R$-trivial, and we finished the proof
of \ref{E7case}.
\end{ntt}

\section{Other homogeneous varieties}
In this section using the results of \cite{Me} and \cite{Ti66}
we finish the proof of the theorem of the introduction.
We start with the following
\begin{lem}\label{ratcorr}
Let $X$ and $Y$ denote projective homogeneous varieties over
a field $k$. Assume $X$ is isotropic over the function field of $Y$
and $Y$ is isotropic over the function field of $X$.
Then the groups of zero-cycles of $X$ and $Y$ are isomorphic.
\end{lem}

\begin{proof}
The fact that $X$ is isotropic over $k(Y)$ is equivalent
to the existence of a rational map $Y\dasharrow X$. Hence,
we have two composable rational maps $f\colon Y\dasharrow X$ and 
$g\colon X\dasharrow Y$,
and the compositions $f\circ g$ and $g\circ f$ correspond to taking
a $k(X)$-point on $X$ and a $k(Y)$-point on $Y$ respectively.

Consider the category of rational correspondences $\RatCor(k)$ introduced
in \cite{Me}. The objects of this category are smooth
projective varieties over $k$ and morphisms $\Mor(X,Y)=\CH_0(Y_{k(X)})$.
The key property of this category is that 
the $\CH_0$-functor factors through it.
Namely, $\CH_0$ is a composition of two functors:
the first is given by taking a graph of a rational map
(any rational map gives rise to a morphism in $\RatCor(k)$),
the second is the realization functor (see \cite[Thm.~3.2]{Me}).

The maps $f$ and $g$ give rise
to the morphisms $[f]$ and $[g]$ in $\RatCor(k)$.
By definition the compositions $[f\circ g]$ and $[g\circ f]$
give the identity maps in the category $\RatCor(k)$.
Hence, the realizations $[f]_*$ and $[g]_*$ give the respective
mutually inverse isomorphisms between $\CH_0(X)$ and $\CH_0(Y)$.
\end{proof}

The next lemma finishes the proof of the theorem of introduction.

\begin{lem}
Let $X$ be an anisotropic projective $G$-homogeneous variety, where
$G$ is a group of type $\F$ or $^1{\E_6}$ with trivial Tits algebras. 
Then $\CHO(X)=0$.
\end{lem}

\begin{proof}
According to Proposition~\ref{propPrime} it is enough to prove
the lemma over fields $k_p$, where $p=2$ or $3$.

Assume $p=2$ and $k=k_2$. 
Consider a Jordan algebra $J$ corresponding to the group $G$.
The algebra $J$ is reduced and, hence,
comes from an octonion algebra $\mathbb{O}$. 
Consider the variety $Y$ of norm zero 
elements of $\mathbb{O}$. It can be easily deduced from Tits diagrams that
the algebras $J_{k(X)}$ and $J_{k(Y)}$ split.
This means the varieties $X_{k(Y)}$ and $Y_{k(X)}$ are isotropic.
By Lemma \ref{ratcorr} $\CHO(X)=\CHO(Y)=0$, where the last equality
holds by \cite{Sw89}.

Assume $p=3$ and $k=k_3$. In this case apply similar arguments
to the pair $X$ and $Y$, where $Y=\OP(J)$.
We obtain $\CHO(X)=\CHO(Y)=0$, where the last equality holds
by Theorem~\ref{cayl}.
\end{proof}

\subsection*{Acknowledgements}
We express our thanks to V.~Chernousov, N.~Karpenko, I.~Panin and M.~Rost
for useful discussions.
We are very grateful to S.~Garibaldi for comments concerning
exceptional groups.

The first and the second authors are sincerely grateful 
to A.~Bak for his 
hospitality, advice, and reading of the manuscript. 
The first author is very thankful to N.~Vavilov for substantial remarks and
would like to acknowledge DAAD {\it Jahresstipendium} and INTAS 03-51-3251. 
The second author appreciates the support of INTAS 00-566.
The last author is grateful to AvH Stiftung and Max-Planck-Institute f\"ur
Mathematik for support.

\bibliographystyle{chicago}

\end{document}